\documentclass[aos,MSNbibl,citesort,dvips]{arximspdf}
\usepackage{graphicx}

\doi{10.1214/11-AOS957}
\volume{40}
\issue{1}
\pubyear{2012}
\firstpage{238}
\lastpage{261}

\makeatletter

\newtheorem{lem}{Lemma}[section]
\newtheorem{prop}[lem]{Proposition}
\newtheorem{theorem}[lem]{Theorem}
\newtheorem{cor}[lem]{Corollary}
\newtheorem{conj}[lem]{Conjecture}

\newproclaim{rem}[lem]{Remark}
\newproclaim{ex}[lem]{Example}

\newcommand{\R}{\mathbb{R}}

\newcommand{\s}{\mathbb{S}}

\makeatother

\begin{document}
\begin{frontmatter}

\title{Geometry of maximum likelihood estimation in~Gaussian graphical models}
\runtitle{Geometry of ML estimation in Gaussian graphical models}

\begin{aug}
\author[A]{\fnms{Caroline} \snm{Uhler}\corref{}\ead[label=e1]{caroline.uhler@ist.ac.at}\ead[label=u1,url]{http://www.carolineuhler.com}}
\runauthor{C. Uhler}
\affiliation{Institute of Science and Technology Austria}
\address[A]{Institute of Science\\
\quad and Technology Austria \\
Am Campus 1\\
A-3400 Klosterneuburg\\
Austria\\
\printead{e1}\\
\printead{u1}} 
\end{aug}

\received{\smonth{4} \syear{2011}}
\revised{\smonth{11} \syear{2011}}

%
\begin{abstract}
We study maximum likelihood estimation in Gaussian graphical models
from a~geometric point of view. An algebraic elimination criterion
allows us to find exact lower bounds on the number of observations
needed to ensure that the maximum likelihood estimator (MLE) exists
with probability one. This is applied to bipartite graphs, grids and
colored graphs. We also study the ML degree, and we present the first
instance of a~graph for which the MLE exists with probability one, even
when the number of observations equals the treewidth.
\end{abstract}

%

\end{frontmatter}

\section{Introduction}

In current statistical applications, we are often faced with problems
involving a~large number of random variables, but only a~small number
of observations (e.g.,~\cite{Hastie}, Chapter 18). This problem arises,
for example, when studying genetic networks: We seek a~model
potentially involving a~vast number of genes, while we are only given
gene expression data of a~few individuals. Gaussian graphical models
have frequently been used to study gene association networks. The
maximum likelihood estimator (MLE) of the covariance matrix is computed
to describe the interaction between different genes (e.g.,~\cite{biology1,biology2}). So the following question is of great interest
from an applied as well as a~theoretical point of view: What is the
minimum number of observations needed to guarantee the existence of the
MLE in a~Gaussian graphical model? It is well known that the MLE exists
with probability one if the number of observations is at least as large
as the number of variables. In this paper we examine the case of fewer
observations.

Gaussian graphical models have been introduced by Dempster~\cite{Dempster} under the name of covariance selection models. Subsequently,
the graphical representation of these models increased in importance.
Lauritzen~\cite{Lauritzen} and Whittaker~\cite{Whittaker} give
introductions to graphical models in general and discuss the connection
between graph and probability distribution for Gaussian graphical models.

Gaussian graphical models are regular exponential families. The
statistical theory of exponential families, as presented, for example,
by Brown~\cite{Brown} or Barndorff-Nielsen~\cite{Barndorff}, is a~strong tool to establish existence and uniqueness of the MLE. The MLE
exists and is unique if and only if the sufficient statistic lies in
the interior of its convex support. We will give a~geometric
description of the convex support of the sufficient statistics and
discuss the connection to the number of samples.

This paper is organized as follows. In Section~\ref{prerequisites}, we
explain the connection between maximum likelihood estimation in
Gaussian graphical models and positive definite matrix completion
problems. In Section~\ref{geometry}, we give a~geometric description of
the problem, and we develop an exact algebraic algorithm to determine
lower bounds on the number of observations needed to ensure existence
of the MLE with probability one. In Section~\ref{secbipartite}, we
discuss the existence of the MLE for bipartite graphs. Section~\ref{smallgraphs} deals with small graphs. The $3\times3$ grid motivated
this paper and is the original problem posed by Steffen Lauritzen
during his lecture on the existence of the MLE in Gaussian graphical
models at the ``Durham Symposium on Mathematical Aspects of Graphical
Models'' on July 8, 2008. The $3\times3$ grid is also the first example
of a~graph for which the MLE exists with probability one even when the
number of observations equals the treewidth of the underlying graph. We
conclude this paper with a~characterization of Gaussian models on
colored 4-cycles in Section~\ref{colored}.

\section{Positive definite matrix completion}
\label{prerequisites}

Let $G=([m],E)$ be an undirected graph on the vertex set $[m]=\{1,\ldots
,m\}$ with edge set $E$. To simplify notation, we assume that $E$
contains all self-loops, that is, $(i,i)\in E$ for all $i\in[m]$. Let
$q$ denote the maximal clique size of $G$. A~graph $G$ is
\textit{chordal} if it contains no chordless cycle of length greater
than 3.
For a~nonchordal graph $G=([m],E)$ one can define a~\textit{chordal
cover} $G^+=([m],E^+)$, which is a~chordal graph satisfying $E\subset
E^+$. We denote its maximal clique size by $q^+$. It is useful to
introduce the notion of a~\textit{minimal chordal cover} $G^*=([m],E^*)$,
where minimality refers to the maximal clique size in the chordal
cover, that is, $q^*=\min(q^+)$. The \textit{treewidth} of a~graph $\tau
(G)$ is defined as
\[
\tau(G)=q^*-1.
\]

A~random vector $X$ taking values in $\mathbb{R}^m$ is said to satisfy
the Gaussian graphical model with graph $G$ if $X$ follows a~multivariate normal distribution obeying the undirected pairwise Markov
property (e.g.,~\cite{Lauritzen,Whittaker}). Assuming the mean to be
zero, this property is as follows:
%
%
\begin{equation}\label{eq1}\qquad
X\sim\mathcal{N}(0,\Sigma),\qquad \Sigma \mbox{ positive definite
with } (\Sigma^{-1})_{ij}=0\quad \forall(i,j)\notin E.
\end{equation}

The results in this paper are based on the assumption that the mean is
a~\textit{known} vector. In particular, we study the case where the mean
is zero. The case where the mean is unknown or partially known is more
complex, since mean and covariance matrix can generally not be
estimated independently. Gehrmann and Lauritzen~\cite{Gehrmann}
describe symmetry relations on the underlying graph which ensure
estimability of the mean vector independently from the true covariance
matrix $\Sigma$.

We denote by $\s^m$ the set of symmetric $m\times m$ matrices and by
$\s ^m_{\succ0}$ the open convex cone of positive definite matrices.
For a~matrix $M\in\s^m$ let $M_G$ denote the \textit{$G$-partial
matrix} consisting of all entries of $M$ corresponding to edges in the
graph~$G$, that is,
\[
M_G=\bigl(M_{ij}\mid(i,j)\in E\bigr).
\]
In particular, all diagonal entries of the partial matrix $M_G$ are
specified, because we assume that the edge set $E$ contains all
self-loops. Equivalently,~$M_G$ is the projection of $M$ onto the
(coordinates indexed by the) edge set of the graph $G$:
\[
\pi_G\dvtx\s^m \rightarrow\R^E,\qquad M \mapsto M_G.
\]

Let $X_1,\ldots,X_n$ denote $n$ independent draws from the distribution
$\mathcal{N}(0,\Sigma)$. Then the \textit{sample covariance matrix} is
given by
\[
S=\frac{1}{n}\sum_{i=1}^n X_i X_i^T.
\]
The $G$-partial sample covariance matrix $S_G$ plays an important role
when studying the existence of the MLE, as seen in the following
theorem first proven by Dempster~\cite{Dempster}.
%
%
\begin{theorem}
\label{Dempster}
In the Gaussian graphical model on $G$, the MLE of the covariance
matrix $\Sigma$ exists if and only if the $G$-partial sample covariance
matrix $S_G$ can be completed to a~positive definite matrix. Then the
MLE $ \hat{\Sigma}$ is the unique completion satisfying $(\hat{\Sigma
}^{-1})_{ij}=0$ for all $(i,j)\notin E$.
\end{theorem}

So checking existence of the MLE in a~Gaussian graphical model is a~special matrix completion problem with a~rank constraint on the partial
matrix given by the number of observations. Matrix completion problems
have been extensively studied, and the following result from~\cite{Grone} is very useful in this context.
%

%
\begin{figure}[b]

\includegraphics{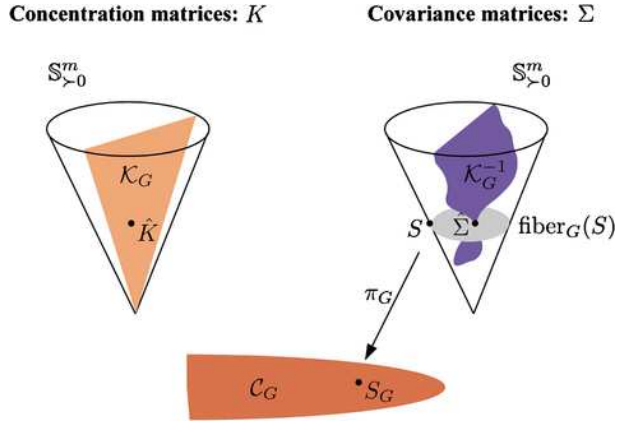}

\caption{Geometry of maximum likelihood estimation in
Gaussian graphical models. The cone $\mathcal{K}_G$ consists of all
concentration matrices in the model, and $\mathcal{K}_G^{-1}$ is the
corresponding set of covariance matrices. The cone of sufficient
statistics $\mathcal{C}_G$ is defined as the projection of $\s^m_{\succ
0}$ onto the edge set of~$G$. It is dual to $\mathcal{K}_G$. Given a~sample covariance matrix $S$, $\mathrm{fiber}_\mathcal{G}(S)$ consists
of all positive definite completions of the $G$-partial matrix $S_G$,
and it intersects $\mathcal{K}_G^{-1}$ in at most one point, namely the
MLE~$\hat{\Sigma}$.} \label{figcones}
\end{figure}

%
\begin{theorem}
\label{thmchordal}
For a~graph $G$ the following statements are equivalent:
\begin{longlist}
\item A~$G$-partial matrix $M_G\in\mathbb{R}^E$ has a~positive
definite completion if and only if all submatrices corresponding to
maximal cliques in $M_G$ are positive definite.
\item $G$ is chordal.\vadjust{\goodbreak}
\end{longlist}
\end{theorem}

By combining Theorems~\ref{Dempster} and~\ref{thmchordal} we
get the following result about the existence of the MLE in Gaussian
graphical models (see also~\cite{Buhl}).
%
%
\begin{cor}
\label{corsamplesize}
If $n\geq q^*$, the MLE exists with probability 1. If $n< q$, the MLE
does not exist.
\end{cor}

Note that chordal graphs have $q^*=q$. Therefore, existence of the MLE
only depends on the number of observations. For nonchordal graphs,
however, there is a~gap $q\leq n < q^*$, in which existence of the MLE
is not well understood. Cycles and wheels (cycles with one additional
completely connected vertex) are the only nonchordal graphs, which
have been studied~\cite{Barrett2,Barrett1,Buhl}. We will extend the
results on cycles and wheels to bipartite graphs $K_{2,m}$ and small
grids.\vspace*{-3pt}

\section{Geometry of maximum likelihood estimation in Gaussian
graphical models}
\label{geometry}

Every concentration matrix (i.e., inverse of a~covariance matrix) in a~Gaussian graphical model satisfies the undirected pairwise Markov
property (\ref{eq1}). The set of all concentration matrices in the
model is a~convex cone
\[
\mathcal{K}_G := \{K\in\s^m_{\succ0} \mid K_{ij} =0,
\forall(i, j) \notin E\}.
\nonumber
\]
Note again that the edge set contains all self-loops, that
is, $(i,i)\in E$ for all \mbox{$i\in[m]$}. By taking the inverse of every
matrix in $\mathcal{K}_G$, we get the set of all covariance matrices in
the model denoted by $\mathcal{K}_G^{-1}$. This is an algebraic variety
intersected with the positive definite cone $\s^m_{\succ0}$ and shown
in purple in Figure~\ref{figcones}.\vadjust{\goodbreak}

In a~Gaussian graphical model, the $G$-partial matrix $S_G$ is a~\textit{minimal sufficient statistic} of a~sample covariance matrix $S$ (e.g.,
\cite{Lauritzen,Whittaker}). So Theorem~\ref{Dempster} has the
following geometric interpretation also explained in Figure~\ref{figcones}:

%
\begin{cor} \label{corfiber}
The MLEs $\hat{\Sigma}$ and $\hat{K}$ exist for a~given sample
covariance matrix $S$ if and only if
\[
\mathrm{fiber}_\mathcal{G}(S) :=
\{ \Sigma\in\s^m_{\succ0} \mid
\Sigma_G = S_G\}
\]
is nonempty, in which case $\mathrm{fiber}_\mathcal{G}(S)$ intersects
$\mathcal{K}_G^{-1}$ in exactly one point, namely the MLE $\hat{\Sigma}$.
\end{cor}

So the MLE $\hat{\Sigma}$ has an algebraic description in terms of the
sufficient statistic~$S_G$, that is, $\hat{\Sigma}$ can be represented
as a~solution to polynomial equations in the sufficient statistic
$S_G$. The maximal degree of these polynomials is called the ML degree.
The ML degree describes the map taking a~sample covariance matrix $S$
to its maximum likelihood estimate $\hat{\Sigma}$ and is studied in
more detail in Section~\ref{secbipartite}.

Applying Corollary~\ref{corfiber}, we can describe the set of all
sufficient statistics for which the MLE exists. We denote this set by
$\mathcal{C}_G$. It is given by the projection of the positive definite
cone $\mathbb{S}^m_{\succ0}$ onto the edge set of the graph~$G$:
\[
\mathcal{C}_G := \pi_G( \s^m_{\succ0} ).
\]
So $\mathcal{C}_G$ is also a~convex cone and shown in dark orange in
Figure~\ref{figcones}. Moreover, we proved in~\cite{ourpaper},
Proposition 2.1, that the cone of sufficient statistics $\mathcal{C}_G$
is the convex dual to the cone of concentration matrices $\mathcal{K}_G$.
%
%
\begin{ex}
\label{excoloredK}
For small-dimensional problems we are able to give a~graphical
representation of the cone of sufficient statistics $\mathcal{C}_G$.
For example, consider the Gaussian graphical model on the bipartite
graph $K_{2,3}$ with concentration matrices of the form
\[
K = \pmatrix{
\lambda_1 & 0 & \lambda_2 & \lambda_3 & \lambda_4\cr
0 & \lambda_1 & \lambda_4 & \lambda_2 & \lambda_3\cr
\lambda_2 & \lambda_4 & \lambda_1 & 0 & 0 \cr
\lambda_3 & \lambda_2 & 0 & \lambda_1 & 0 \cr
\lambda_4 & \lambda_3 & 0 & 0 & \lambda_1},\qquad
\raisebox{-27pt}{
\includegraphics{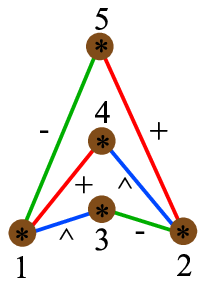}
}
\]
Note that in order to reduce the number of parameters and be able to
draw~$\mathcal{C}_G $ in three-dimensional space, we assume additional
equality constraints on the nonzero entries of the concentration
matrix, represented by the graph coloring above. Such colored Gaussian
graphical models, where the coloring represents equality constraints on
the concentration matrix, are called RCON-models and have been
introduced in~\cite{Hojsgaard}.

%
\begin{figure}
\begin{tabular}{c c c c c c}

\includegraphics{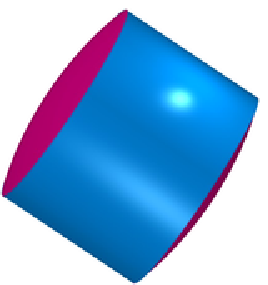}
 & & \includegraphics{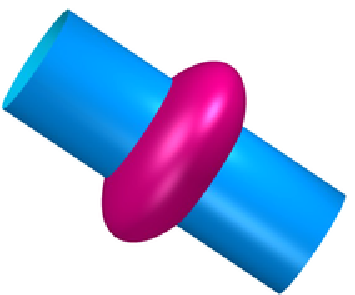} \\
(a) & & (b)\\[6pt]

\includegraphics{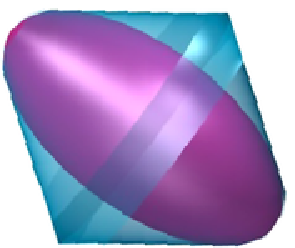}
 & \includegraphics{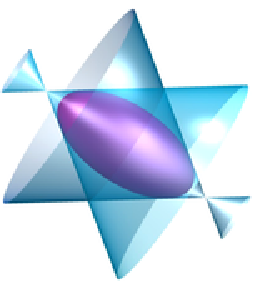} & \includegraphics{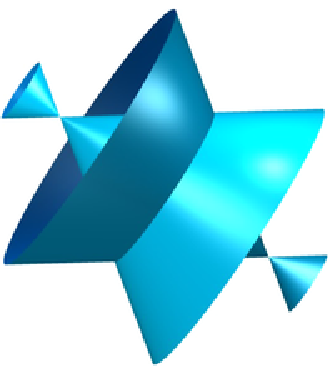} \\
(c) & (d) & (e)
\end{tabular}
\caption{These pictures illustrate the convex geometry of maximum
likelihood estimation for Gaussian graphical models. The cone of
concentration matrices $ \mathcal{K}_G $ is shown in \textup{(a)}, its
algebraic boundary in \textup{(b)}, the dual cone of sufficient statistics in
\textup{(c)} and its algebraic boundary in \textup{(d)} and \textup{(e)},
where \textup{(d)} is the transparent version of \textup{(e)}.} \label{KG}
\end{figure}

Without loss of generality we can rescale $K$ and assume that all
diagonal entries are one. The cone of concentration matrices\vadjust{\goodbreak} $
\mathcal{K}_G $ for this model is shown in Figure~\ref{KG}(a). Its
algebraic boundary is described by $\{\operatorname{det}(K) = 0\}$ and
is shown in Figure~\ref{KG}(b). In this example, the determinant
factors into two components, a~cylinder and an ellipsoid. Dualizing the
boundary of $\mathcal{K}_G $ by the algorithm described in our previous
paper (\cite{ourpaper}, Proposition 2.4) results in the hypersurface
shown in Figure~\ref{KG}(e). The double cone is dual to the cylinder in
Figure~\ref{KG}(b). By making the double cone transparent as shown in
Figure~\ref{KG}(d), we see the enclosed ellipsoid, which is dual to the
ellipsoid in Figure~\ref{KG}(b). The cone of sufficient statistics
$\mathcal{C}_G$ is shown in Figure~\ref{KG}(c). The MLE exists if and
only if the sufficient statistic lies in the interior of this convex
body. Using the elimination criterion of Theorem~\ref{thmelimination},
we can show that the MLE exists with probability one already for one
observation.
\end{ex}

In this paper, we examine the existence of the MLE for $n$ observations
in the range $q\leq n < q^*$, for which the existence of the MLE is not
well understood. Geometrically, we look at the manifold of rank $n$
matrices on the boundary of the cone $\s^m_{\succeq0}$. In general,
its projection
%
%
\begin{equation}
\label{projection}
\pi_G\bigl(\{M\in\s^m_{\succeq0} \mid\mathrm{rk}(M)=n \}
\bigr)
\end{equation}
lies in the topological closure of the cone $\mathcal{C}_G$. The MLE
exists with probability one for $n$ observations if and only if the
projection (\ref{projection}) lies in the interior of $\mathcal{C}_G$.

Based on the geometric interpretation of maximum likelihood estimation
in Gaussian graphical models, we can derive a~sufficient condition for
the existence of the MLE. The following algebraic elimination criterion
can be used as an algorithm to establish existence of the MLE with
probability one for $n$ observation.
%
%
\begin{theorem}[(Elimination criterion)]
\label{thmelimination}
Let $I_{G,n}$ be the elimination ideal obtained from the ideal of
$(n+1) \times(n+1)$-minors of a~symmetric $m\times m$ matrix $S$ of unknowns by eliminating all unknowns
corresponding to nonedges of the graph~$G$. If $I_{G,n}$ is the zero ideal,
then the MLE exists with probability one for $n$ observations.
\end{theorem}
\begin{pf}
The variety corresponding to the ideal of $(n+1) \times(n+1)$-minors
of a~symmetric $m\times m$ matrix $S$ of unknowns consists of all $m\times
m$ matrices of rank at most $n$. Eliminating all unknowns
corresponding to nonedges of the graph $G$ results in the elimination
ideal $I_{G,n}$ (see, e.g.,~\cite{Cox}) and is geometrically equivalent
to a~projection onto the cone of sufficient statistics $\mathcal{C}_G$.
Let $V$ be the variety corresponding to the elimination ideal~$I_{G,n}$.
We denote by $k$ its dimension and by $\mu$ a~$k$-dimensional Lebesgue measure. The MLE exists with probability one
for $n$ observations if
\[
\mu(V\cap\partial\mathcal{C}_G )=0,
\]
where $\partial\mathcal{C}_G$ denotes the boundary of the cone of
sufficient statistics $\mathcal{C}_G$.

If $I_{G,n}$ is the zero ideal, then the variety $V$ is
full-dimensional, and its dimension $\dim(V)=k=\dim(\mathcal{C}_G)$. So
if we assume that $\mu(V\cap\partial\mathcal{C}_G )>0$, then $\mu(
\partial\mathcal{C}_G )>0$, which is a~contradiction to $\dim(
\partial\mathcal{C}_G )<k$.
\end{pf}

For small examples, the elimination ideal $I_{G,n}$ can be computed, for
example, using \texttt{Macaulay2}~\cite{Macaulay2}, a~software system for
research in algebraic geometry. If $I_{G,n}$ is not the zero ideal,
then an analysis of polynomial inequalities is required. One needs to
carefully examine how the components of $V$ are located. The argument
is subtle because the algebraic boundary of $\mathcal{C}_G$ may in fact
intersect the interior of $\mathcal{C}_G$. So even if the projection
$V$ is a~component of the algebraic boundary of $\mathcal{C}_G$, the
MLE might still exist with positive probability. We will encounter and
describe such an example in detail in Section~\ref{colored}.

\section{Bipartite graphs}
\label{secbipartite}

\begin{figure}

\includegraphics{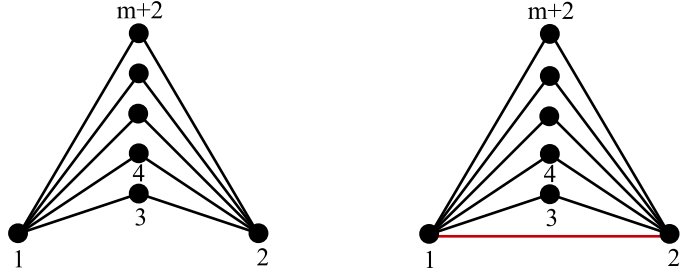}

\caption{Bipartite graph $K_{2,m}$ (left) and minimal chordal cover of
$K_{2,m}$ (right).}
\label{bipartite}
\end{figure}

In this section, we first derive the MLE existence results for
bipartite graphs $K_{2,m}$ paralleling the results on cycles proven by
Buhl~\cite{Buhl}. Let the graph $K_{2,m}$ be labeled as shown in Figure
\ref{bipartite}. A~minimal chordal cover is given in Figure
\ref{bipartite} (right). As for cycles, for bipartite graphs~$K_{2,m}$
we have $q=2$ and $q^*=3$. Therefore only the case of $n=2$
observations is interesting.

%

%
\begin{figure}[b]

\includegraphics{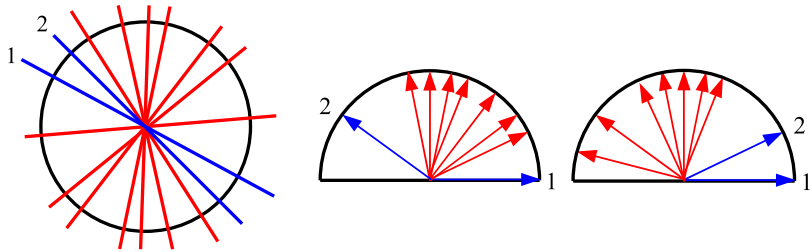}

\caption{The MLE on $K_{2,m}$ exists in the following situations. Lines
and data vectors corresponding to the variables 1 and 2 are drawn in
blue. Lines and data vectors corresponding to the variables $3,4,\ldots
,m+2$ are drawn in red.}
\label{MLEK2m}
\end{figure}

Let $X_1$ and $X_2$ denote two independent samples from the
distribution $\mathcal{N}_{m+2}(0$, $\Sigma)$, which obeys the undirected\vadjust{\goodbreak}
pairwise Markov property on $K_{2,m}$. We denote by $X$ the
$(m+2)\times2$ data matrix consisting of the two samples~$X_1$ and
$X_2$ as columns. The rows of $X$ are denoted by $x_1,\ldots, x_{m+2}$.
Similarly as for cycles in~\cite{Buhl}, we will describe a~criterion on
the configuration of data vectors $x_1,\ldots, x_{m+2}$ ensuring the
existence of the MLE. Our proof is essentially the same argument as
used by Buhl~\cite{Buhl} for cycles. The following characterization of
positive definite matrices of size $3\times3$ proven in~\cite{Barrett1} will be helpful in this context.
%
%
\begin{lem}
\label{lemarccos}
The matrix
\[
\pmatrix{
1 & \cos(\alpha) & \cos(\beta)\cr
\cos(\alpha) & 1 & \cos(\gamma)\cr
\cos(\beta) & \cos(\gamma) & 1}
\]
with $0<\alpha,\beta,\gamma<\pi$ is positive definite if and
only if
\[
\alpha<\beta+\gamma,\qquad \beta<\alpha+\gamma,\qquad \gamma<\alpha
+\beta,\qquad \alpha+\beta+\gamma<2\pi.
\]
\end{lem}
%
%
\begin{prop}
\label{propbipartite}
The MLE on the graph $K_{2,m}$ exists with probability one for $n\geq
3$ observations, and the MLE does not exist for $n<2$ observations. For
$n=2$ observations the MLE exists if and only if the lines generated by
$x_1$ and $x_2$ are direct neighbors [see Figure~\ref{MLEK2m} (left)].
\end{prop}
\begin{pf}
Because the problem of existence of the MLE is a~positive definite
matrix completion problem, we can rescale and rotate the data vectors
$x_1,\ldots, x_{m+2}$ (i.e., perform an orthogonal transformation)
without changing the problem. So without loss of generality we can
assume that the vectors $x_1,\ldots, x_{m+2}\in\R^2$ have length one,
lie in the upper unit half circle and $x_1=(1,0)$. We need to prove
that the MLE exists if and only if the data configuration is as shown
in Figure~\ref{MLEK2m} (middle) or (right).

Let $\theta_{ij}$ denote the angle between vector $x_i$ and $x_j$. Then
the $K_{2,m}$-partial sample covariance matrix $S_{K_{2,m}}$ is of the form
\[
\left(
\begin{tabular}{cc|cccc}
1 & $\star$ & $\cos(\theta_{13})$ & $\cos(\theta_{14})$ & $\cdots$ &
$\cos(\theta_{1,m+2})$ \\
$\star$ & 1 & $\cos(\theta_{23})$ & $\cos(\theta_{24})$ & $\cdots$ &
$\cos(\theta_{2,m+2})$ \\ \hline
$\cos(\theta_{13})$ & $\cos(\theta_{23})$ & 1 & $\star$ & $\cdots$ &
$\star$\\
$\cos(\theta_{14})$ & $\cos(\theta_{24})$ & $\star$ & 1 & $\ddots$&
$\vdots$\\
$\vdots$ & $\vdots$ & $\vdots$ & $\ddots$&$\ddots$& $\star$\\
$\cos(\theta_{1,m+2})$ & $\cos(\theta_{2,m+2})$ & $\star$ & $\cdots
$&$\star$ & 1\\
\end{tabular}
\right).
\]
We put stars ($\star$) at all positions not corresponding to edges in
the graph. The stars represent the entries of the sample covariance
matrix which are not part of the sufficient statistics.

The graph $K_{2,m}$ can be extended to a~chordal graph by adding one
edge as shown in Figure~\ref{bipartite} (right). So by Theorem
\ref{thmchordal}, $S_{K_{2,m}}$ can be extended to a~positive definite
matrix if and only if the $(1,2)$ entry of $S_{K_{2,m}}$ can be
completed in such a~way that all the submatrices corresponding to
maximal cliques are positive definite. This is equivalent to the
existence of $\rho\in\R$ with $0<\rho<\pi$ such that
\[
\pmatrix{
1 & \cos(\rho) & \cos(\theta_{1i}) \cr
\cos(\rho) & 1 & \cos(\theta_{2i}) \cr
\cos(\theta_{1i}) & \cos(\theta_{2i}) & 1}
\succ0 \qquad\mbox{for all } i\in\{3, 4,\ldots,m+2\}.
\]
By Lemma~\ref{lemarccos} this occurs if and only if
\[
\left.\matrix{
\theta_{1i}-\theta_{2i}\cr
\theta_{2i} -\theta_{1i}}\right\} <\rho<
\left\{\matrix{
\theta_{1j}+\theta_{2j}\cr 2\pi-\theta_{1j} -\theta_{2j}}\right.
\qquad\mbox{for all } i,j\in\{3,4,\ldots,m+2\},
\]
which is equivalent to
%
%
\begin{equation}
\label{eqbipartite}
2\theta_{ai}<\theta_{1i}+\theta_{2i}+\theta_{1j}+\theta_{2j}<2\pi
+2\theta_{ai}
\end{equation}
for all $a\in\{1,2\}$, $i,j\in\{3,4,\ldots,m+2\}$. We distinguish two cases.
\begin{longlist}[\textit{Case} 2.]
\item[\textit{Case} 1.] There is a~vector $x_j$ lying between $x_1$ and $x_2$,
which implies that $\theta_{1j}+\theta_{2j}=\theta_{12}$. If there was
a~vector $x_i$, $i\neq j$, which does not lie between $x_1$ and $x_2$, then
\[
\theta_{1j}+\theta_{2j}+\theta_{1i}+\theta_{2i}=2\theta_{1i},
\]
which is a~contradiction to (\ref{eqbipartite}). Hence all vectors
$x_3, x_4, \ldots x_{m+2}$ lie between~$x_1$ and $x_2$, in which case
\[
\theta_{1i}+\theta_{2i}+\theta_{1j}+\theta_{2j}=2\theta_{12},
\]
and inequality (\ref{eqbipartite}) is satisfied.
\item[\textit{Case} 2.] The vectors $x_1$ and $x_2$ are direct neighbors, which
implies that $\theta_{1i}+\theta_{2i}=\theta_{12}+2\theta_{2i}$ for all
$i\in\{3, 4, \ldots,m+2\}$, in which case inequality (\ref
{eqbipartite}) is satisfied.
\end{longlist}

This proves that for two observations, the MLE exists if and
only if the data configuration is as shown in Figure~\ref{MLEK2m}
(middle) or (right).
\end{pf}

The geometric explanation of what is happening in this example is that
the projection of the positive definite matrices of rank 2 intersects
the interior and the boundary of the cone of sufficient statistics
$\mathcal{C}_G$ with positive measure. The sufficient statistics
originating from data vectors, where lines~1 and 2 are neighbors, lie
in the interior of $\mathcal{C}_G$. If lines 1 and 2 are not neighbors,
the corresponding sufficient statistics lie on the boundary of the
cone~$\mathcal{C}_G$, and the MLE does not exist. A~similar situation is
encountered in Example~\ref{exFret} and depicted in Figure~\ref{figfrets}.

It is worth remarking that if the $m+2$ variables are independent, we
can compute the probability of existence of the MLE by a~combinatorial
argument. In this case, the probability that the MLE exists is given by
\[
\frac{2m!}{(m+1)!} = \frac{2}{m+1}.
\]

A~different approach to gaining a~better understanding of maximum
likelihood estimation in Gaussian graphical models is to study the ML
degree of the underlying graph. The map taking a~sample covariance
matrix $S$ to its maximum likelihood estimate $\hat\Sigma$ is an
algebraic function, and its degree is the ML degree of the model. See
\cite{Oberwolfach}, Definition 2.1.4. The ML degree represents the
algebraic complexity of the problem of finding the MLE. This suggests
that a~larger ML degree results in a~more difficult MLE existence
problem. We proved in~\cite{ourpaper} that the ML degree is one if and
only if the underlying graph is chordal. It is conjectured in~\cite{Oberwolfach},
Section 7.4, that the ML degree of the cycle grows
exponentially in the cycle length. An interesting contrast to the cycle
conjecture is the following theorem, where we prove that the ML degree
for bipartite graphs $K_{2,m}$ grows linearly in the number of variables.
%
%
\begin{theorem}
\label{MLdegreebipartite}
$\!\!\!$In a~Gaussian graphical model with underlying graph~$K_{2,m}$ the ML
degree is $2m+1$.
\end{theorem}
\begin{pf}
Given a~generic matrix $S\in\mathbb{S}^{m+2}$, we fix $\Sigma\in\mathbb
{S}^{m+2}$ with entries $\Sigma_{ij}=S_{ij}$ for $(i,j)\in E$ and
unknowns $\Sigma_{12}=\Sigma_{21}=y$ and $\Sigma_{ij}=z_{ij}$ for all
other $(i,j)\notin E$. We denote by $K=\Sigma^{-1}$ the corresponding
concentration\vadjust{\goodbreak} matrix. The ML degree of $K_{2,m}$ is the number of
complex solutions to 
\[
(\Sigma^{-1})_{ij}=0 \qquad\mbox{for all } (i,j)\notin E.
\]

Let $A$ denote the set consisting of the two distinguished vertices $\{
1,2\}$, and let $B=V\setminus A$. In the following we will use the
block structure
\[
\Sigma=\pmatrix{ \Sigma_{AA} & \Sigma_{AB} \cr
\Sigma_{BA} & \Sigma_{BB}},\qquad
K=\pmatrix{ K_{AA} & K_{AB} \cr K_{BA} & K_{BB}}.
\]
For example, for the graph $K_{2,5}$ the corresponding covariance
matrix $\Sigma$ and concentration matrix $K$ are of the form
\begin{eqnarray*}
\Sigma &=& \left(
\begin{tabular}{cc|ccc}
1 & $y$ & $S_{13}$ & $S_{14}$ & $S_{15}$ \\
$y$ & $1$ & $S_{23}$ & $S_{24}$ & $S_{25}$ \\ \hline
$S_{13}$ & $S_{23}$ & 1 & $z_{34}$ & $z_{35}$ \\
$S_{14}$ & $S_{24}$ & $z_{34}$ & 1 & $z_{45}$ \\
$S_{15}$ & $S_{25}$ & $z_{35}$ & $z_{45}$ & 1
\end{tabular}
\right) ,\\
K &=& \left(
\begin{tabular}{cc|ccc}
$K_{11}$ & $0$ & $K_{13}$ & $K_{14}$ & $K_{15}$ \\
0 & $K_{22}$ & $K_{23}$ & $K_{24}$ & $K_{25}$ \\ \hline
$K_{13}$ & $K_{23}$ & $K_{33}$ & 0 & 0 \\
$K_{14}$ & $K_{24}$ & 0 & $K_{44}$ & 0 \\
$K_{15}$ & $K_{25}$ & 0 & 0 & $K_{55}$
\end{tabular}
\right).
\end{eqnarray*}
Note that the block $K_{BB}$ is a~diagonal matrix. Hence the Schur complement
\[
\Sigma_{BB}-\Sigma_{BA}\Sigma_{AA}^{-1}\Sigma_{AB}
\]
is also a~diagonal matrix. Writing out the off-diagonal entries of this
matrix results in the following expression for the variables $z$ in
terms of the variable~$y$:
\[
z_{ij}=-\frac{1}{1-y^2}\bigl(y(S_{1i}S_{2j}+S_{1j}S_{2i})
-S_{1i}S_{1j}-S_{2i}S_{2j}\bigr).
\]
Setting the minor $M_{12}$ of $\Sigma$ to zero results in the last
equation of the form
%
%
\begin{equation}
\label{eqpolybipartite}
y\det(\Sigma_{BB})+(\mbox{polynomial in $z$ of degree $m-1$})=0.
\end{equation}
We note that $\det(\Sigma_{BB})$ is a~polynomial in $z$ of degree $m$,
where the degree 0 term is~1. So by multiplying equation (\ref
{eqpolybipartite}) with $(1-y^2)^m$, we get a~degree $2m+1$ equation
in y and therefore $2m+1$ complex solutions for $y$. For each solution
of $y$ we get one solution for the variables $z$, which proves that the
ML degree of $K_{2,m}$ is $2m+1$.
\end{pf}

Bipartite graphs and cycles are classes of graphs with $q=2$ and
$q^*=3$. What can we say about such graphs in general regarding the
existence of the MLE for two observations? A~related question has been
studied from a~purely algebraic point of\vadjust{\goodbreak} view in~\cite{Barrett2}. A~cycle-completable graph is defined to be a~graph such that every
partial matrix $M_G$ has a~positive definite completion if and only if
$M_G$ is positive definite on all submatrices corresponding to maximal
cliques in the graph, and all submatrices corresponding to cycles in
the graph can be completed to a~positive definite matrix. It is shown
in~\cite{Barrett2} that a~graph is cycle-completable if and only if
there is a~chordal cover with no new 4-clique.

Buhl~\cite{Buhl} studied cycles from a~more statistical point of view
and described a~criterion on the data vectors for the existence of the
MLE for two observations. Combining the results of~\cite{Barrett2} and
\cite{Buhl}, we get the following result:

%
\begin{cor}
\label{corq3}
Let $G$ be a~graph with $q=2$ and $q^*\geq3$. Then the following
statements are equivalent:
\begin{longlist}
\item For $n=2$ observations, the MLE exists if and only if
Buhl's cycle condition is satisfied on every induced cycle.
\item$q^*=3$.
\end{longlist}
%
\end{cor}

This result solves the problem of existence of the MLE for all graphs
with $q=2$ and $q^*=3$. Note that Corollary~\ref{corq3} is more
general than Proposition~\ref{propbipartite}. The proof, however, is
more involved and less constructive.

For bipartite graphs $K_{3,m}$ the situation is more complicated and we
do not yet have results similar to Proposition~\ref{propbipartite} and
Theorem~\ref{MLdegreebipartite}. We will nevertheless describe some
preliminary results.

%
%
\begin{figure}

\includegraphics{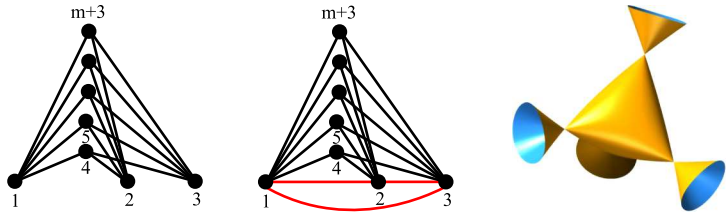}

\caption{Bipartite graph $K_{3,m}$ (left) and minimal chordal cover of
$K_{3,m}$ (middle). The tetrahedron-shaped pillow consisting of all
positive semidefinite $3\times3$ matrices with ones on the diagonal is
shown in the right figure.}
\label{bipartite3}
\end{figure}

Let the graph $K_{3,m}$ be labeled as shown in Figure
\ref{bipartite3}. A~minimal chordal cover is given in Figure
\ref{bipartite3} (middle). Hence, $q=2$ and $q^*=4$. The convex body
shown in Figure~\ref{bipartite3} (right) consists of all positive
semidefinite $3\times3$ matrices with ones on the diagonal. We call it
the tetrahedron-shaped pillow. We will prove that the existence of the
MLE is equivalent to a~nonempty intersection of such inflated and
shifted tetrahedron-shaped pillows.
%
%
\begin{cor}\label{cor4.5}
The MLE on the graph $K_{3,m}$ exists if and only if the $m$ inflated
and shifted tetrahedron-shaped pillows corresponding to the maximal
cliques in a~minimal chordal cover of $K_{3,m}$ shown in Figure
\ref{bipartite3} (middle) have nonempty intersections.\vadjust{\goodbreak}
\end{cor}
\begin{pf}
Applying Theorem~\ref{thmchordal} in a~similar way as in the proof of
Theorem~\ref{propbipartite}, the partial covariance matrix
$S_{K_{3,m}}$ can be extended to a~positive definite matrix if and only
if the entries corresponding to the missing edges $(1,2)$, $(1,3)$ and
$(2,3)$ can be completed in such a~way that all the submatrices
corresponding to maximal cliques in the minimal chordal cover (Figure
\ref{bipartite3}, middle) are positive definite. This is equivalent to
the existence of $x, y, z\in\R$ with \mbox{$-1<x,y,z<1$} such that
%
%
\begin{equation}\label{eq4matrix}
\pmatrix{
1 & s_{1i} & s_{2i} & s_{3i}\cr
s_{1i} & 1 & x & y\cr
s_{2i} & x & 1 & z\cr
s_{3i} & y & z & 1}
\succ0 \qquad\mbox{for all } i\in\{4, 5,\ldots,m+3\},
\end{equation}
where $s_{ai}$, $a\in\{1,2,3\}$, $i\in\{4,5,\ldots, m+3\}$ are the
sufficient statistics corresponding to edges in the bipartite graph
$K_{3,m}$. Using Schur complements and rescaling, (\ref{eq4matrix})
holds if and only if
%
%
\begin{equation}
\label{eq5matrix}
\pmatrix{
1 & x_i & y_i\cr
x_i & 1 & z_i\cr
y_i & z_i & 1}
\succ0 \qquad\mbox{for all } i\in\{4, 5,\ldots,m+3\},
\end{equation}
where
\begin{eqnarray*}
x_i&=&\frac{x-s_{1i}s_{2i}}{\sqrt{1-s_{1i}^2}\sqrt{1-s_{2i}^2}},\qquad
y_i=\frac{y-s_{1i}s_{3i}}{\sqrt{1-s_{1i}^2}\sqrt{1-s_{3i}^2}},\\
z_i&=&\frac{x-s_{2i}s_{3i}}{\sqrt{1-s_{2i}^2}\sqrt{1-s_{3i}^2}}.
\end{eqnarray*}
So the MLE exists if and only if the inflated and shifted
tetrahedron-shaped pillows corresponding to the inequalities in (\ref
{eq5matrix}) have nonempty intersection.
\end{pf}

We used the software package \texttt{Macaulay2} to compute the ML degree
of~$K_{3,m}$ for $m\leq4$. It is an open problem to find a~general
formula or a~recurrence relation for the ML degree of $K_{l,m}$, where
$l\geq3$.\vspace*{6pt}

\begin{center}
\begin{tabular}{c|cccc}
$m$ & 1 & 2 & 3 & 4 \\ \hline
ML degree & 1 & 7 & 57 & 131
\end{tabular}
\end{center}

\section{Small graphs}
\label{smallgraphs}

In this section we analyze the $3\times3$ grid in particular and
complete the discussion of~\cite{ourpaper} with the number of
observations and the corresponding existence probability of the MLE for
all graphs with 5 or less vertices.

%
\begin{figure}

\includegraphics{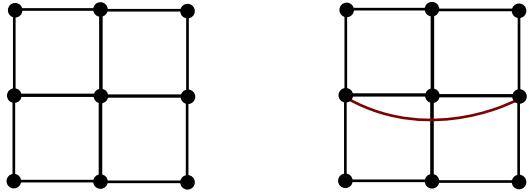}

\caption{$3\times3$ grid $\mathcal{H}$ (left) and grid with additional
edge $\mathcal{H}'$ (right).} \label{figgrid}
\end{figure}

%
\begin{figure}[b]

\includegraphics{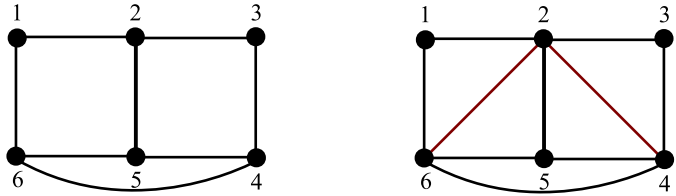}

\caption{Graph $\mathcal{G}$ (left) and minimal chordal
cover of $\mathcal{G}$ (right).} \label{ourgraph}
\end{figure}

The $3\times3$ grid is shown in Figure~\ref{figgrid} (left) and has
$q=2$ and $q^*=4$. This example represents the starting point of this
paper and is the original problem posed by Steffen Lauritzen during his
lecture at the ``Durham Symposium on Mathematical Aspects of Graphical
Models'' in 2008. As a~preparation, we first discuss the existence of
the MLE for the graph $\mathcal{G}$ on six vertices shown in Figure~\ref{ourgraph}. The graph $\mathcal{G}$ also has $q=2$ and $q^*=4$, and is
the first example for which we can prove that the bound $n\geq q^*$ for
the existence of the MLE with probability one is not tight and that the
MLE can exist
with probability one, even when the number of observations equals the treewidth.

%
\begin{theorem}\label{thmourgraph}
The MLE on the graph $\mathcal{G}$ (Figure~\ref{ourgraph}, left)
exists with probability one for $n= 3$ observations.
\end{theorem}
\begin{pf}
We compute the ideal $I_{\mathcal{G},3}$ by eliminating the variables
$s_{13}, s_{15},\allowbreak s_{16}$, $s_{24}, s_{26}, s_{34}, s_{35}$ from the
ideal of $4\times4$ minors of the matrix $S$ given in (\ref
{matrixxS}). This results in the zero ideal, which by Theorem~\ref{thmelimination} completes the proof.
\end{pf}
%
%
\begin{rem}\label{remalg}
Theorem~\ref{thmourgraph} is equivalent to the following purely
algebraic statement. Let
%
%
\begin{equation}
\label{matrixxS}
S=\pmatrix{ 1 & s_{12} & s_{13} & s_{14} & s_{15} & s_{16} \cr
s_{12} & 1 & s_{23} & s_{24} & s_{25} & s_{26} \cr
s_{13} & s_{23} & 1 & s_{34} & s_{35} & s_{36} \cr
s_{14} & s_{24} & s_{34} & 1 & s_{45} & s_{46} \cr
s_{15} & s_{25} & s_{35} & s_{45} & 1 & s_{56} \cr
s_{16} & s_{26} & s_{36} & s_{46} & s_{56} & 1 }
\in\mathbb{S}^6_{\succeq0}
\end{equation}
with $\operatorname{rank}(S)=3$. Then there exist $x, y, a, b, c, d, e
\in\R
$ such that
\[
S'=\pmatrix{ 1 & s_{12} & a& s_{14} & b & c \cr
s_{12} & 1 & s_{23} & x & s_{25} & y \cr
a & s_{23} & 1 & d & e & s_{36} \cr
s_{14} & x & d & 1 & s_{45} & s_{46} \cr
b & s_{25} & e & s_{45} & 1 & s_{56} \cr
c & y & s_{36} & s_{46} & s_{56} & 1}
\in\mathbb{S}^6_{\succ0}.
\]
So any partial matrix of rank 3 with specified entries at all positions
corresponding to edges in $\mathcal{G}$ can be completed to a~positive
definite matrix.
\end{rem}
%
%
\begin{cor}
\label{cor3by3}
Let $\mathcal{H}$ be the $3\times3$-grid shown in Figure~\ref{figgrid}. Then the MLE on $\mathcal{H}$ exists with probability one
for $n\geq3$ observations, and the MLE does not exist for $n<2$ observations.
\end{cor}
\begin{pf}
First note that Groebner bases computations are extremely memory
intensive and the elimination ideal $I_{\mathcal{H},3}$ cannot be
computed directly due to insufficient memory. We solve this problem by
gluing together smaller graphs. The probability of existence of the MLE
for the $3\times3$ grid $\mathcal{H}$ is at least as large as the
existence probability when the underlying graph is $\mathcal{H}'$. The
graph $\mathcal{H}'$ is a~clique sum of two graphs of the form $\mathcal
{G}$, for which the MLE existence probability is one for $n\geq3$.
\end{pf}

This example shows that although we are not able to compute the
elimination ideal for large graphs directly, the algebraic elimination
criterion (Theorem~\ref{thmelimination}) is still useful also in this
situation. We can study small graphs with the elimination criterion and
glue them together using clique sums to build larger graphs.

For two observations on the $3\times3$ grid, the cycle conditions are
necessary but not sufficient for the existence of the MLE (Corollary
\ref{corq3}). Unlike for bipartite graphs $K_{2,m}$, the existence of
the MLE does not only depend on the ordering of the lines corresponding
to the data vectors in $\mathbb{R}^2$. By simulations with the
\texttt{Matlab} software \texttt{cvx}~\cite{CVX}, one can easily find orderings
for which the MLE sometimes exists and sometimes does not. Finding a~necessary and sufficient criterion for the existence of the MLE for two
observations remains an open problem.

%
\begin{table}
\tablewidth=280pt
\caption{This table shows the number of observations (obs.) and the
corresponding MLE existence probability for all nonchordal graphs on~5
or fewer vertices}
\label{tablesmallgraphs}
\begin{tabular*}{\tablewidth}{@{\extracolsep{\fill}}lcccc@{}}
\hline
\textbf{Graph} $\bolds{G}$ & \textbf{1 obs.} & \textbf{2 obs.}
& \textbf{3 obs.} & $\bolds{\geq}$\textbf{4 obs.}\\
\hline
(a)
\raisebox{-10pt}{\includegraphics{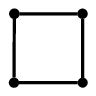}}
 & No & $p\in(0,1)$ & $p=1$ & $p=1$ \\[20pt]
(b)
\raisebox{-10pt}{\includegraphics{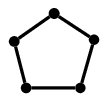}}
 & No & $p\in(0,1)$ & $p=1$ & $p=1$ \\[20pt]
(c)
\raisebox{-10pt}{\includegraphics{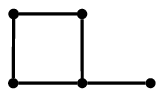}}
 & No & $p\in(0,1)$ & $p=1$ & $p=1$ \\[20pt]
(d)
\raisebox{-10pt}{\includegraphics{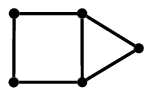}}
 & No & No & $p=1$ & $p=1$ \\[20pt]
(e)
\raisebox{-10pt}{\includegraphics{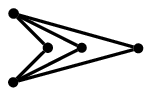}}
 & No & $p\in(0,1)$ & $p=1$ & $p=1$ \\[20pt]
(f)
\raisebox{-10pt}{\includegraphics{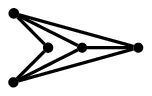}}
 & No & No & $p=1$ & $p=1$ \\[20pt]
(g)
\raisebox{-10pt}{\includegraphics{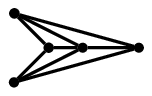}}
 & No & No & $p\in(0,1)$ & $p=1$\\
\hline
\end{tabular*}
\end{table}

We now complete the discussion of~\cite{ourpaper} with the number of
observations and the corresponding existence probability of the MLE for
all graphs with~5 or less vertices. All nonchordal graphs with 5 or
less vertices are shown in Table~\ref{tablesmallgraphs}. The 4-cycle
and 5-cycle in (a) and (b) are covered by Buhl's results~\cite{Buhl}. The
graphs in (c) and (d) are clique sums of two graphs and therefore
completable if and only if the submatrices corresponding to the two
subgraphs are completable. Graph (e) is the bipartite graph $K_{2,3}$
and covered by Theorem~\ref{propbipartite}.\vadjust{\goodbreak} For the graph in (f) $q=3$
and $q^*=4$. Applying the elimination criterion from Theorem~\ref{thmelimination} shows that three observations are sufficient for the
existence of the MLE. The last example, the 5-wheel in graph (g),
is also covered by Buhl's results~\cite{Buhl}.

\section{Colored Gaussian graphical models}
\label{colored}

For some applications, symmetries in the underlying Gaussian graphical
model can be assumed. Adding symmetry to the conditional
independence restrictions of a~graphical model reduces the number of
parameters and in some cases also the number of observations needed for
the existence of the MLE. The symmetry restrictions can be represented
by a~graph coloring, where edges, or vertices, respectively, have the
same coloring if the corresponding elements of the concentration matrix
are equal. Such models are called RCON-models~\cite{Hojsgaard}. We
discussed such a~model earlier in Example~\ref{excoloredK}.

We denote the uncolored graph by $G$ and the colored graph by $\mathcal
{G}$. Note that in this section the graph $G$ does not contain any
self-loops. Let the vertices be colored with $p$
different colors and the edges with $q$ different colors:
\begin{eqnarray*}
V & = & V_1\sqcup V_2\sqcup\cdots\sqcup V_p,\qquad p\leq|V|, \\
E & = & E_1\sqcup E_2\sqcup\cdots\sqcup E_q,\qquad q\leq|E|.
\end{eqnarray*}
Then the set of all concentration matrices $\mathcal{K}_{\mathcal{G}}$
consists of all positive definite matrices satisfying:
\begin{itemize}
\item$ K_{\alpha\beta} = 0 $ for any pair of vertices $\alpha, \beta
$ that do not
form an edge in $G$.
\item
$ K_{\alpha\alpha}=K_{\beta\beta}$ for any pair of vertices $\alpha,
\beta$ in a~common vertex color class~$V_i$.
\item
$ K_{\alpha\beta}=K_{\gamma\delta}$ for any pair of edges $(\alpha
,\beta), (\gamma,\delta)$
in a~common edge color class~$E_j$.
\end{itemize}
This means that also for RCON-models the set $\mathcal{K}_{\mathcal
{G}}$ is defined by linear equations on the concentration matrix $K$.
So the geometry of maximum likelihood estimation is the same as that
explained in Section~\ref{geometry}, and it is straightforward to
derive the equivalent of Theorem~\ref{Dempster} for colored Gaussian
graphical models.
%
%
\begin{theorem} \label{Dempstercolored}
In a~colored Gaussian graphical model on $\mathcal{G}$ the MLE of the
covariance matrix $\Sigma$ exists if and only if there is a~positive
definite matrix $\tilde{\Sigma}$ such that
\[
\sum_{\alpha\in V_i} \tilde{\Sigma}_{\alpha\alpha}=\sum_{\alpha\in
V_i} S_{\alpha\alpha} \quad\mbox{and}\quad \sum_{(\alpha, \beta
)\in E_j} \tilde{\Sigma}_{\alpha\beta}=\sum_{(\alpha, \beta)\in E_j}
S_{\alpha\beta}
\]
for all vertex color classes $V_1, \ldots, V_p$ and all edge color
classes $E_1, \ldots, E_q$. Then the MLE $ \hat{\Sigma}$ is the unique
completion with $ (\hat{\Sigma}^{-1})_{\alpha\alpha}=(\hat{\Sigma
}^{-1})_{\beta\beta}$ for any pair of vertices $\alpha, \beta$ in a~common vertex color class $V_i$, $ (\hat{\Sigma}^{-1})_{\alpha\beta
}=(\hat{\Sigma}^{-1})_{\gamma\delta}$ for any pair of edges $(\alpha
,\beta), (\gamma,\delta)$ in a~common edge color class $E_j$, and
$(\hat{\Sigma}^{-1})_{\alpha\beta}=0$ for all $(\alpha, \beta)\notin E$.
\end{theorem}
%
\begin{ex}[(Frets's heads)]
\label{exFret}
We revisit the heredity study of head dimensions known as \textit{Frets's
heads} reported in~\cite{Frets}. Part of the original data are the
length and breadth of the heads of 25 pairs of first and second sons.
This data set was also discussed in~\cite{MKB,ourpaper}. The data
supports the following colored Gaussian graphical model, where the
joint distribution remains the same when the two sons are exchanged:
\[
K = \pmatrix{ \lambda_1 & \lambda_3 & 0 & \lambda_4 \cr
\lambda_3 & \lambda_1 & \lambda_4 & 0\cr
0 & \lambda_4 & \lambda_2 & \lambda_5 \cr
\lambda_4 & 0 & \lambda_5 & \lambda_2 },\qquad
\raisebox{-19pt}{
\includegraphics{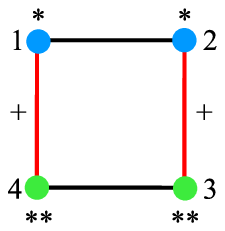}
}
\]
In this graph, variable 1 corresponds to the length of the first son's
head, variable 2 to the length of the second son's head,\vadjust{\goodbreak} variable 3 to
the breadth of the second son's head and variable 4 to the breadth of
the first son's head. Color classes consisting only of one edge (or
vertex) are displayed in black.

%
\begin{figure}

\includegraphics{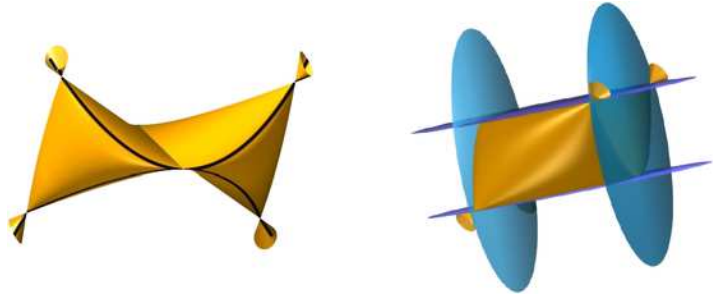}

\caption{All possible sufficient statistics from one observation are
shown on the left. The cone of sufficient statistics is shown on the
right.} \label{figfrets}
\end{figure}

Given a~sample covariance matrix $S=(s_{ij})$, the
five sufficient statistics for this model according to the graph
coloring are
\begin{eqnarray*}
t_1 &=& s_{11}+s_{22},\qquad t_2=s_{33}+s_{44},\qquad
t_3=2s_{12},\\
t_4 &=& 2(s_{23}+s_{14}),\qquad t_5=2s_{34}.
\end{eqnarray*}

The algebraic boundary of the cone of sufficient statistics
$\mathcal{C}_{\mathcal{G}}$ is computed in~\cite{ourpaper} and given by
the polynomial
\begin{eqnarray*}
H_{\mathcal{G}}&=&
(t_1-t_3)\cdot(t_1+t_3)\cdot(t_2-t_5)\cdot(t_2+t_5)\\
&&{}\times
(4t_2^2 t_3^2-4 t_1 t_2
t_4^2+t_4^4+8t_1t_2t_3t_5-4t_3t_4^2t_5+4t_1^2t_5^2).
\end{eqnarray*}
For two observations the elimination ideal $I_{\mathcal{G},2}$ is the
zero ideal. Therefore, the MLE exists with probability 1 for two or
more observations in this model. For one observation we get
\[
I_{\mathcal{G},1}=\langle4t_2^2 t_3^2-4 t_1 t_2
t_4^2+t_4^4+8t_1t_2t_3t_5-4t_3t_4^2t_5+4t_1^2t_5^2\rangle,
\]
which corresponds to one of the components of the algebraic boundary of
the cone of sufficient statistics. In this example, the algebraic
boundary of the cone of sufficient statistics intersects its interior.
This is illustrated in Figure~\ref{figfrets}. In order to get a~graphical representation in three-dimensional space, we fixed $t_3$ and
$t_5$. The variety corresponding to $I_{\mathcal{G},1}$ is shown on the
left. We call this hypersurface the bow tie. The cone of sufficient
statistics~$\mathcal{C}_{\mathcal{G}}$ is the convex hull of the bow
tie and shown in Figure~\ref{figfrets} (right). Its boundary consists
of four planes corresponding to the components $t_1-t_3$, $t_1+t_3$,
$t_2-t_5$ and $t_2+t_5$ shown in blue, and the bows of the bow tie
shown in yellow. The black curves show where the planes touch the bow
tie. Note that the upper and lower two triangles of the bow tie lie in
the interior of~$\mathcal{C}_{\mathcal{G}}$. Only the two bows are part
of the boundary of $\mathcal{C}_{\mathcal{G}}$. So the MLE exists if
the sufficient statistic lies\vadjust{\goodbreak} on one of the triangles of the bow tie,
and it does not exist if the sufficient statistic lies on one of the
bows of the bow tie. Consequently, for one observation the MLE exists
with probability strictly between 0 and~1.

A~different approach is to run simulations, for example, using
\texttt{cvx}. We can generate vectors of length four and compute the
MLE by
solving a~convex optimization problem. If \texttt{cvx} finds a~solution,
the MLE exists. For this example, however, \texttt{cvx} sometimes does
not find a~solution, which supports the hypothesis that the MLE exists
with probability strictly between 0 and 1 for one observation.
In the following, we give a~formal proof by characterizing the set of
vectors in $\mathbb{R}^4$ for which the MLE exists/does not exist.


For this example, we can exactly characterize not just the sufficient
statistics, but also the observations, for which the MLE exists. In
other words, we can characterize the observations whose sufficient
statistics lie on the triangles of the bow tie. First, note that by
exchanging variables 1 and 2 and simultaneously exchanging variables 3
and 4, we get the same model. This means that from one observation
$X_1=(x_1,x_2,x_3,x_4)$ we can generate a~second observation
$X_2=(x_2,x_1,x_4,x_3)$. So the resulting data matrix is given
by\looseness=-1
\[
X=\pmatrix{ x_1&x_2\cr x_2&x_1\cr x_3&x_4\cr x_4&x_3}.
\]\looseness=0
Applying Buhl's result about two observations on a~Gaussian cycle~\cite{Buhl}, the MLE exists if and only if the lines corresponding to the vectors
\[
y_1=\pmatrix{x_1\cr x_2},\qquad y_2=\pmatrix{x_2\cr x_1},\qquad
y_3=\pmatrix{x_3\cr
x_4},\qquad y_4=\pmatrix{x_4\cr x_3}
\]
are not graph consecutive. This is the case if and only if
%
%
\begin{equation}
\label{ineqfrets}
|x_1|>|x_2| \mbox{ and } |x_3|>|x_4| \quad\mbox{or}\quad
|x_1|<|x_2| \mbox{ and } |x_3|<|x_4|.
\end{equation}
Hence, the MLE for one observation exists if and only if the data is
inconsistent, meaning that the head of the first (second) son is longer
than the head of the second (first) son, but the breadth is smaller. In
this situation the corresponding sufficient statistics lie on the
triangles of the bow tie in Figure~\ref{figfrets}. Otherwise the
corresponding sufficient statistics lie on the bows of the bow tie. If
$\Sigma$ is diagonal, the MLE exists with probability 0.5, since all
configurations in (\ref{ineqfrets}) have the same probability.
\end{ex}

In our previous paper~\cite{ourpaper} we found the defining polynomial
$\mathcal{H}_{\mathcal{G}}$ of the cone of sufficient statistics for
all colored Gaussian graphical models on the 4-cycle, which have the
property that edges in the same color class connect the same vertex
color classes. Such models have been studied in~\cite{Hojsgaard} and
are of special interest, because they are invariant
under rescaling of variables in the same vertex color class. In Tables
\ref{RCOR1} and~\ref{RCOP},\vadjust{\goodbreak} we complete the discussion of
\cite{ourpaper} with the number of observations and the corresponding
existence probability of the MLE.

%
\begin{table}
\tablewidth=280pt
\caption{Results on the number of observations and the MLE existence
probability for all colored
Gaussian graphical models with some symmetry restrictions
(namely, edges in the same color
class connect the same vertex color classes) on the 4-cycle}
\label{RCOR1}
\begin{tabular*}{\tablewidth}{@{\extracolsep{\fill}}lcccc@{}}
\hline
\textbf{Graph} & $\bolds{K}$ & \textbf{1 obs.} & \textbf{2 obs.}
& $\bolds{\geq}$\textbf{3 obs.} \\
\hline
(1)
\raisebox{-12pt}{\includegraphics{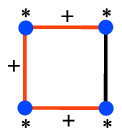}}
 & $
\pmatrix{ \lambda_1 & \lambda_2 & 0 & \lambda_2\cr
\lambda_2 & \lambda_1 & \lambda_3 & 0\cr
0 & \lambda_3 & \lambda_1 & \lambda_2\cr
\lambda_2 & 0 & \lambda_2 & \lambda_1}$ & $p=1$ & $p=1$ & $p=1$\\ [34pt]
(2)
\raisebox{-12pt}{\includegraphics{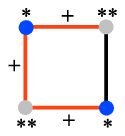}}
 & $
\pmatrix{ \lambda_1 & \lambda_3 & 0 & \lambda_3\cr
\lambda_3 & \lambda_2 & \lambda_4 & 0\cr
0 & \lambda_4 & \lambda_1 & \lambda_3\cr
\lambda_3 & 0 & \lambda_3 & \lambda_2 }$ & $p=1$ & $p=1$ & $p=1$\\[34pt]
(3)
\raisebox{-12pt}{\includegraphics{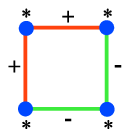}}
 & $
\pmatrix{ \lambda_1 & \lambda_2 & 0 & \lambda_2\cr
\lambda_2 & \lambda_1 & \lambda_3 & 0\cr
0 & \lambda_3 & \lambda_1 & \lambda_3\cr
\lambda_2 & 0 & \lambda_3 & \lambda_1}$ & $p=1$ & $p=1$ & $p=1$\\[34pt]
(4)
\raisebox{-12pt}{\includegraphics{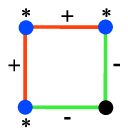}}
 & $
\pmatrix{ \lambda_1 & \lambda_3 & 0 & \lambda_3\cr
\lambda_3 & \lambda_1 & \lambda_4 & 0\cr
0 & \lambda_4 & \lambda_2 & \lambda_4\cr
\lambda_3 & 0 & \lambda_4 & \lambda_1 }$ & $p=1$ & $p=1$ & $p=1$\\[34pt]
(5)
\raisebox{-12pt}{\includegraphics{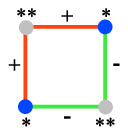}}
 & $
\pmatrix{ \lambda_1 & \lambda_3 & 0 & \lambda_3\cr
\lambda_3 & \lambda_2 & \lambda_4 & 0\cr
0 & \lambda_4 & \lambda_1 & \lambda_4\cr
\lambda_3 & 0 & \lambda_4 & \lambda_2 }$ & $p=1$ & $p=1$ & $p=1$\\[34pt]
(6)
\raisebox{-12pt}{\includegraphics{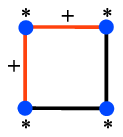}}
 & $
\pmatrix{ \lambda_1 & \lambda_2 & 0 & \lambda_2\cr
\lambda_2 & \lambda_1 & \lambda_3 & 0\cr
0 & \lambda_3 & \lambda_1 & \lambda_4\cr
\lambda_2 & 0 & \lambda_4 & \lambda_1}$ & $p=1$ & $p=1$ & $p=1$\\[34pt]
(7)
\raisebox{-12pt}{\includegraphics{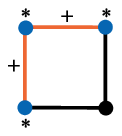}}
 & $
\pmatrix{ \lambda_1 & \lambda_3 & 0 & \lambda_3\cr
\lambda_3 & \lambda_1 & \lambda_4 & 0\cr
0 & \lambda_4 & \lambda_2 & \lambda_5\cr
\lambda_3 & 0 & \lambda_5 & \lambda_1}$ & $p\in(0,1)$ & $p=1$ & $p=1$\\[34pt]
(8)
\raisebox{-12pt}{\includegraphics{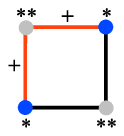}}
 & $
\pmatrix{ \lambda_1 & \lambda_3 & 0 & \lambda_3\cr
\lambda_3 & \lambda_2 & \lambda_4 & 0\cr
0 & \lambda_4 & \lambda_1 & \lambda_5\cr
\lambda_3 & 0 & \lambda_5 & \lambda_2}$ & $p\in(0,1)$ & $p=1$ & $p=1$\\[34pt]
(9)
\raisebox{-12pt}{\includegraphics{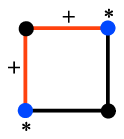}}
 & $
\pmatrix{ \lambda_1 & \lambda_4 & 0 & \lambda_4\cr
\lambda_4 & \lambda_2 & \lambda_5 & 0\cr
0 & \lambda_5 & \lambda_3 & \lambda_6\cr
\lambda_4 & 0 & \lambda_6 & \lambda_2}$ & No? & $p=1$ & $p=1$\\
\hline
\end{tabular*}
\end{table}

\setcounter{table}{1}
\begin{table}
\tablewidth=280pt
\caption{(Continued)}
\begin{tabular*}{\tablewidth}{@{\extracolsep{\fill}}lcccc@{}}
\hline
\textbf{Graph} & $\bolds{K}$ & \textbf{1 obs.} & \textbf{2 obs.}
& $\bolds{\geq}$\textbf{3 obs.} \\
\hline
(10)
\raisebox{-12pt}{\includegraphics{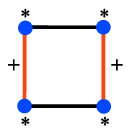}}
 & $
\pmatrix{\lambda_1 & \lambda_2 & 0 & \lambda_3\cr
\lambda_2 & \lambda_1 & \lambda_3 & 0\cr
0 & \lambda_3 & \lambda_1 & \lambda_4\cr
\lambda_3 & 0 & \lambda_4 & \lambda_1}$ & $p=1$ & $p=1$ & $p=1$\\[34pt]
(11)
\raisebox{-12pt}{\includegraphics{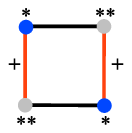}}
 & $
\pmatrix{ \lambda_1 & \lambda_3 & 0 & \lambda_4\cr
\lambda_3 & \lambda_2 & \lambda_4 & 0\cr
0 & \lambda_4 & \lambda_1 & \lambda_5\cr
\lambda_4 & 0 & \lambda_5 & \lambda_2}$ & No? & $p=1$ & $p=1$\\[34pt]
(12)
\raisebox{-12pt}{\includegraphics{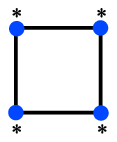}}
 & $\pmatrix{ \lambda_1 &
\lambda_2 & 0 &
\lambda_5\cr
\lambda_2 & \lambda_1 & \lambda_3 & 0\cr
0 & \lambda_3 & \lambda_1 & \lambda_4\cr
\lambda_5 & 0 & \lambda_4 & \lambda_1 }$ & $p=1$ & $p=1$ & $p=1$\\[34pt]
%
(13)
\raisebox{-12pt}{\includegraphics{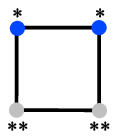}}
 & $\pmatrix{ \lambda_1&
\lambda_3& 0&
\lambda_6\cr
\lambda_3& \lambda_1& \lambda_4& 0\cr
0& \lambda_4& \lambda_2& \lambda_5\cr
\lambda_6& 0& \lambda_5& \lambda_2}$ & $p\in(0,1)$ & $p=1$ & $p=1$\\[34pt]
(14)
\raisebox{-12pt}{\includegraphics{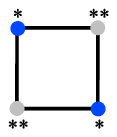}}
 & $\pmatrix{ \lambda_1 &
\lambda_3 & 0 &
\lambda_6\cr
\lambda_3 & \lambda_2 & \lambda_4 & 0\cr
0 & \lambda_4 & \lambda_1 & \lambda_5\cr
\lambda_6 & 0 & \lambda_5 & \lambda_2}$ & No? & $p=1$ & $p=1$\\[34pt]
(15)
\raisebox{-12pt}{\includegraphics{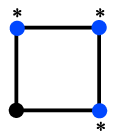}}
 & $\pmatrix{ \lambda_1 &
\lambda_3 & 0 &
\lambda_6\cr
\lambda_3 &\lambda_1 &\lambda_4 &0\cr
0 &\lambda_4 &\lambda_1 &\lambda_5\cr
\lambda_6 &0 &\lambda_5 &\lambda_2}$ & $p\in(0,1)$ & $p=1$ & $p=1$\\[34pt]
(16)
\raisebox{-12pt}{\includegraphics{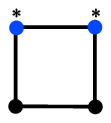}}
 & $\pmatrix{ \lambda_1
&\lambda_4& 0&
\lambda_7\cr
\lambda_4 &\lambda_1 &\lambda_5 &0\cr
0 &\lambda_5 &\lambda_2 &\lambda_6\cr
\lambda_7 &0& \lambda_6 &\lambda_3}$ & No & $p=1$ & $p=1$\\[34pt]
(17)
\raisebox{-12pt}{\includegraphics{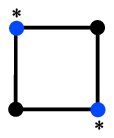}}
 & $\pmatrix{ \lambda
_1&\lambda_4& 0&
\lambda_7\cr
\lambda_4 &\lambda_2 &\lambda_5 &0\cr
0 &\lambda_5 &\lambda_1 &\lambda_6\cr
\lambda_7 &0& \lambda_6 &\lambda_3}$ & No? & $p=1$ & $p=1$\\[34pt]
(18)
\raisebox{-12pt}{\includegraphics{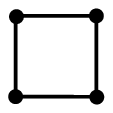}}
 & $\pmatrix{ \lambda_1&
\lambda_5& 0&
\lambda_8\cr
\lambda_5 &\lambda_2 &\lambda_6 &0\cr
0 &\lambda_6 &\lambda_3 &\lambda_7\cr
\lambda_8 &0 &\lambda_7 &\lambda_4}$ & No & $p\in(0,1)$ & $p=1$\\
\hline
\end{tabular*}
\end{table}
%
%
%
\begin{table}
\tablewidth=280pt
\caption{All RCOP-models (introduced in \protect\cite{Hojsgaard}), that is,
graphs with an additional permutation property on the $4$-cycle}
\label{RCOP}
\begin{tabular*}{\tablewidth}{@{\extracolsep{\fill}}lcccc@{}}
\hline
\textbf{Graph} & $\bolds{K}$ & \textbf{1 obs.} & \textbf{2 obs.}
& $\bolds{\geq}$\textbf{3 obs.} \\
\hline
(1)
\raisebox{-12pt}{\includegraphics{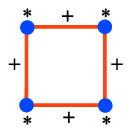}}
 & $\pmatrix{ \lambda_1 &
\lambda_2 & 0 & \lambda_2
\cr
\lambda_2 & \lambda_1 & \lambda_2 & 0 \cr
0 & \lambda_2 & \lambda_1 & \lambda_2 \cr
\lambda_2 & 0 & \lambda_2 & \lambda_1}$ & $p=1$ & $p=1$ & $p=1$\\[34pt]
(2)
\raisebox{-12pt}{\includegraphics{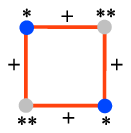}}
 & $\pmatrix{ \lambda_1 &
\lambda_3 & 0 &
\lambda_3\cr
\lambda_3 & \lambda_2 & \lambda_3 & 0\cr
0 & \lambda_3 & \lambda_1 & \lambda_3\cr
\lambda_3 & 0 & \lambda_3 & \lambda_2}$ & $p=1$ & $p=1$ & $p=1$ \\[34pt]
(3)
\raisebox{-12pt}{\includegraphics{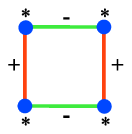}}
 & $\pmatrix{ \lambda_1 &
\lambda_2 & 0 &
\lambda_3\cr
\lambda_2 & \lambda_1 & \lambda_3 & 0\cr
0 & \lambda_3 & \lambda_1 & \lambda_2\cr
\lambda_3 & 0 & \lambda_2 & \lambda_1}$ & $p=1$ & $p=1$ & $p=1$\\[34pt]
(4)
\raisebox{-12pt}{\includegraphics{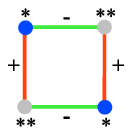}}
 & $\pmatrix{ \lambda_1 &
\lambda_3 & 0 &
\lambda_4\cr
\lambda_3 & \lambda_2 & \lambda_4 & 0\cr
0 & \lambda_4 & \lambda_1 & \lambda_3\cr
\lambda_4 & 0 & \lambda_3 & \lambda_2 }$ & $p=1$ & $p=1$ & $p=1$ \\[34pt]
(5)
\raisebox{-12pt}{\includegraphics{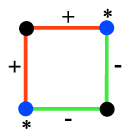}}
 & $\pmatrix{ \lambda_1 &
\lambda_4 & 0 &
\lambda_4\cr
\lambda_4 & \lambda_2 & \lambda_5 & 0\cr
0 & \lambda_5 & \lambda_3 & \lambda_5\cr
\lambda_4 & 0 & \lambda_5 & \lambda_2}$ & $p=1$ & $p=1$ & $p=1$ \\[34pt]
(6)
\raisebox{-12pt}{\includegraphics{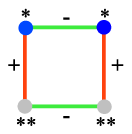}}
 & $\pmatrix{ \lambda_1 &
\lambda_3 & 0 &
\lambda_4\cr
\lambda_3 & \lambda_1 & \lambda_4 & 0\cr
0 & \lambda_4 & \lambda_2 & \lambda_5\cr
\lambda_4 & 0 & \lambda_5 & \lambda_2 }$ & $p\in(0,1)$ & $p=1$ &
$p=1$\\
\hline
\end{tabular*}
\vspace*{-6pt}
\end{table}

For every colored 4-cycle, we computed the elimination ideal
$I_{\mathcal{G},n}$ for $n=1,2,3$. If it is the zero ideal, we know
from Theorem~\ref{thmelimination} that the MLE exists with probability
one. If $I_{\mathcal{G},n}$ is nonzero, we run simulations using
\texttt{cvx}. If we find examples for which the MLE exists and other examples
for which the MLE does not exist, it indicates that the MLE exists with
probability strictly between 0 and 1 for $n$ observations. In cases
where simulations do not yield any counterexamples, we need to prove
that the MLE does indeed not exist by carefully analyzing the
components corresponding to the ideal~$I_{\mathcal{G},n}$. This is the
case for one observation on the graphs (9), (11), (14) and~(17). Note that
the graphical models (9) and (11) are sub-models of (14) and (17). So if we
prove that the MLE does not exist for one observation on the graphs~(9)
and (11), this follows also for the graphs (14) and (17).

If the cone $\mathcal{C}_{\mathcal{G}}$ for the graphs (9) and (11) is a~basic open semialgebraic set (see, e.g.,~\cite{basic}), then $\mathcal
{C}_{\mathcal{G}}$ does not meet its algebraic boundary, and the MLE
does not exist for one observation. So we end with the following
conjecture which would answer the question marks in
Table~\ref{RCOR1}:\vspace*{-2pt}
%
%
\begin{conj}
The cones $\mathcal{C}_{\mathcal{G}}$ corresponding to the graphs (9)
and~(11) are basic open semialgebraic sets.\vspace*{-2pt}
\end{conj}
%

\section{Conclusion}
In this paper, we explained the geometry of maximum likelihood
estimation in Gaussian graphical models. The geometric picture can be
translated into an algebraic criterion (Theorem~\ref{thmelimination}),
which allows us to find exact lower bounds on the number of
observations needed for the existence of the MLE (with probability 1).
Theorem~\ref{thmelimination} holds for any Gaussian graphical model.
However, the practical implementation of Theorem~\ref{thmelimination}
is based on Groebner bases computations, which are extremely memory
intensive. Theorem~\ref{thmourgraph} and Corollary~\ref{cor3by3} show
the power but also the limitations of computational algebraic geometry.
We are, in practice, only able to apply the algebraic elimination
criterion directly to very small graphs. One way of getting results for
larger graphs is to find a~clique decomposition into small subgraphs,
which can be handled individually. A~different future line of research
is to use the small examples to understand the existence of the MLE
asymptotically. If we fix a~class of graphs, for example, cycles or
grids, what can we say about the existence of the MLE as the number of
vertices tends to infinity? Medium-sized graphs, however, remain
untouched by both approaches, and finding the minimum number of
observations needed for the existence of the MLE for such graphs is an
interesting open problem.\vspace*{-2pt}

\section*{Acknowledgments}

I wish to thank Bernd Sturmfels for many helpful discussions and
Steffen Lauritzen for introducing me to the problem of the existence of
the MLE in Gaussian graphical models. I would also like to thank two
referees who provided helpful comments on the original version of this
paper.\vspace*{-2pt}


%

\printaddresses


\begin{thebibliography}{23}

\bibitem{basic}
%
\begin{barticle}[mr]
\bauthor{\bsnm{Acquistapace},~\bfnm{F.}\binits{F.}},
\bauthor{\bsnm{Broglia},~\bfnm{F.}\binits{F.}} \AND
\bauthor{\bsnm{V{\'e}lez},~\bfnm{M.~P.}\binits{M.~P.}}
(\byear{1999}).
\btitle{Basicness of semialgebraic sets}.
\bjournal{Geom. Dedicata}
\bvolume{78}
\bpages{229--240}.
\bid{doi={10.1023/A:1005123421867}, issn={0046-5755}, mr={1725377}}
\bptok{imsref}%
\end{barticle}
%
\endbibitem

\bibitem{Barndorff}
%
\begin{bbook}[mr]
\bauthor{\bsnm{Barndorff-Nielsen},~\bfnm{Ole}\binits{O.}}
(\byear{1978}).
\btitle{Information and Exponential Families in Statistical Theory}.
\bpublisher{Wiley}, \baddress{Chichester}.
\bid{mr={0489333}}
\bptok{imsref}%
\end{bbook}
%
\endbibitem

\bibitem{Barrett1}
%
\begin{barticle}[mr]
\bauthor{\bsnm{Barrett},~\bfnm{Wayne}\binits{W.}},
\bauthor{\bsnm{Johnson},~\bfnm{Charles~R.}\binits{C.~R.}} \AND
\bauthor{\bsnm{Tarazaga},~\bfnm{Pablo}\binits{P.}}
(\byear{1993}).
\btitle{The real positive definite completion problem for a~simple cycle}.
\bjournal{Linear Algebra Appl.}
\bvolume{192}
\bpages{3--31}.
\bid{doi={10.1016/0024-3795(93)90234-F}, issn={0024-3795}, mr={1236734}}
\bptok{imsref}%
\end{barticle}\vadjust{\goodbreak}
%
\endbibitem

\bibitem{Barrett2}
%
\begin{barticle}[mr]
\bauthor{\bsnm{Barrett},~\bfnm{Wayne~W.}\binits{W.~W.}},
\bauthor{\bsnm{Johnson},~\bfnm{Charles~R.}\binits{C.~R.}} \AND
\bauthor{\bsnm{Loewy},~\bfnm{Raphael}\binits{R.}}
(\byear{1996}).
\btitle{The real positive definite completion problem: Cycle completability}.
\bjournal{Mem. Amer. Math. Soc.}
\bvolume{122}
\bpages{viii+69}.
\bid{issn={0065-9266}, mr={1342017}}
\bptok{imsref}%
\end{barticle}
%
\endbibitem

\bibitem{Brown}
%
\begin{bbook}[mr]
\bauthor{\bsnm{Brown},~\bfnm{Lawrence~D.}\binits{L.~D.}}
(\byear{1986}).
\btitle{Fundamentals of Statistical Exponential Families with
Applications in
Statistical Decision Theory}.
\bseries{Institute of Mathematical Statistics Lecture Notes---Monograph Series}
\bvolume{9}.
\bpublisher{IMS}, \baddress{Hayward, CA}.
\bid{mr={0882001}}
\bptok{imsref}%
\end{bbook}
%
\endbibitem

\bibitem{Buhl}
%
\begin{barticle}[mr]
\bauthor{\bsnm{Buhl},~\bfnm{S{\o}ren~L.}\binits{S.~L.}}
(\byear{1993}).
\btitle{On the existence of maximum likelihood estimators for graphical
{G}aussian models}.
\bjournal{Scand. J. Stat.}
\bvolume{20}
\bpages{263--270}.
\bid{issn={0303-6898}, mr={1241392}}
\bptok{imsref}%
\end{barticle}
%
\endbibitem

\bibitem{Cox}
%
\begin{bbook}[mr]
\bauthor{\bsnm{Cox},~\bfnm{David}\binits{D.}},
\bauthor{\bsnm{Little},~\bfnm{John}\binits{J.}} \AND
\bauthor{\bsnm{O'Shea},~\bfnm{Donal}\binits{D.}}
(\byear{1997}).
\btitle{Ideals, Varieties, and Algorithms: An Introduction to Computational Algebraic
Geometry and Commutative Algebra}.
\bpublisher{Springer}, \baddress{New York}.
\bptok{imsref}%
\end{bbook}
%
\endbibitem

\bibitem{Dempster}
%
\begin{barticle}[author]
\bauthor{\bsnm{Dempster},~\bfnm{A.~P.}\binits{A.~P.}}
(\byear{1972}).
\btitle{Covariance selection}.
\bjournal{Biometrics}
\bvolume{28}
\bpages{157--175}.
\bptok{imsref}%
\end{barticle}
%
\endbibitem

\bibitem{Oberwolfach}
%
\begin{bbook}[mr]
\bauthor{\bsnm{Drton},~\bfnm{Mathias}\binits{M.}},
\bauthor{\bsnm{Sturmfels},~\bfnm{Bernd}\binits{B.}} \AND
\bauthor{\bsnm{Sullivant},~\bfnm{Seth}\binits{S.}}
(\byear{2009}).
\btitle{Lectures on Algebraic Statistics}.
\bseries{Oberwolfach Seminars}
\bvolume{39}.
\bpublisher{Birkh\"auser}, \baddress{Basel}.
\bid{doi={10.1007/978-3-7643-8905-5}, mr={2723140}}
\bptok{imsref}%
\end{bbook}
%
\endbibitem

\bibitem{Frets}
%
\begin{barticle}[author]
\bauthor{\bsnm{Frets},~\bfnm{G.~P.}\binits{G.~P.}}
(\byear{1921}).
\btitle{Heredity of head form in man}.
\bjournal{Genetica}
\bvolume{3}
\bpages{193--400}.
\bptok{imsref}%
\end{barticle}
%
\endbibitem

\bibitem{Gehrmann}
%
\begin{bmisc}[author]
\bauthor{\bsnm{Gehrmann},~\bfnm{H.}\binits{H.}} \AND
\bauthor{\bsnm{Lauritzen},~\bfnm{S.~L.}\binits{S.~L.}}
(\byear{2011}).
\bhowpublished{Estimation of means in graphical Gaussian models with symmetries.
Preprint. Available at
\texttt{\href{http://arxiv.org/abs/1101.3709}{http://arxiv.org/abs/}
\href{http://arxiv.org/abs/1101.3709}{1101.3709}}.}
\bptok{imsref}%
\end{bmisc}
%
\endbibitem

\bibitem{CVX}
%
\begin{bmisc}[author]
\bauthor{\bsnm{Grant},~\bfnm{M.}\binits{M.}} \AND
\bauthor{\bsnm{Boyd},~\bfnm{S.}\binits{S.}}
\bhowpublished{CVX, a~Matlab software for disciplined convex programming.
Available at
\url{http://cvxr.com/cvx/}.}
\bptok{imsref}%
\end{bmisc}
%
\endbibitem

\bibitem{Macaulay2}
%
\begin{bmisc}[author]
\bauthor{\bsnm{Grayson},~\bfnm{D.~R.}\binits{D.~R.}} \AND
\bauthor{\bsnm{Stillman},~\bfnm{M.~E.}\binits{M.~E.}}
\bhowpublished{Macaulay2, a~software system for research in algebraic geometry.
Available at \url{http://www.math.uiuc.edu/Macaulay2/}.}
\bptok{imsref}%
\end{bmisc}
%
\endbibitem

\bibitem{Grone}
%
\begin{barticle}[mr]
\bauthor{\bsnm{Grone},~\bfnm{Robert}\binits{R.}},
\bauthor{\bsnm{Johnson},~\bfnm{Charles~R.}\binits{C.~R.}},
\bauthor{\bparticle{de} \bsnm{S{\'a}},~\bfnm{Eduardo~M.}\binits
{E.~M.}} \AND
\bauthor{\bsnm{Wolkowicz},~\bfnm{Henry}\binits{H.}}
(\byear{1984}).
\btitle{Positive definite completions of partial {H}ermitian matrices}.
\bjournal{Linear Algebra Appl.}
\bvolume{58}
\bpages{109--124}.
\bid{doi={10.1016/0024-3795(84)90207-6}, issn={0024-3795}, mr={0739282}}
\bptok{imsref}%
\end{barticle}
%
\endbibitem

\bibitem{Hastie}
%
\begin{bbook}[author]
\bauthor{\bsnm{Hastie},~\bfnm{T.}\binits{T.}},
\bauthor{\bsnm{Tibshirani},~\bfnm{R.}\binits{R.}} \AND
\bauthor{\bsnm{Friedman},~\bfnm{J.}\binits{J.}}
(\byear{2009}).
\btitle{The Elements of Statistical Learning},
\bedition{2nd} ed.
\bseries{Springer Series in Statistics}.
\bpublisher{Springer}, \baddress{New York}.
\bid{mr={2722294}}
\bptok{imsref}%
\end{bbook}
%
\endbibitem

\bibitem{Hojsgaard}
%
\begin{barticle}[mr]
\bauthor{\bsnm{H{\o}jsgaard},~\bfnm{S{\o}ren}\binits{S.}} \AND
\bauthor{\bsnm{Lauritzen},~\bfnm{Steffen~L.}\binits{S.~L.}}
(\byear{2008}).
\btitle{Graphical {G}aussian models with edge and vertex symmetries}.
\bjournal{J. R. Stat. Soc. Ser. B Stat. Methodol.}
\bvolume{70}
\bpages{1005--1027}.
\bid{doi={10.1111/j.1467-9868.2008.00666.x}, issn={1369-7412}, mr={2530327}}
\bptok{imsref}%
\end{barticle}
%
\endbibitem

\bibitem{Lauritzen}
%
\begin{bbook}[mr]
\bauthor{\bsnm{Lauritzen},~\bfnm{Steffen~L.}\binits{S.~L.}}
(\byear{1996}).
\btitle{Graphical Models}.
\bpublisher{Clarendon}, \baddress{Oxford}.
\bid{mr={1419991}}
\bptok{imsref}%
\end{bbook}
%
\endbibitem

\bibitem{MKB}
%
\begin{bbook}[author]
\bauthor{\bsnm{Mardia},~\bfnm{K.~V.}\binits{K.~V.}},
\bauthor{\bsnm{Kent},~\bfnm{J.~T.}\binits{J.~T.}} \AND
\bauthor{\bsnm{Bibby},~\bfnm{J.~M}\binits{J.~M.}}
(\byear{1979}).
\btitle{Multivariate Analysis}.
\bpublisher{Academic Press}, \baddress{London}.
\bptok{imsref}%
\end{bbook}
%
\endbibitem

\bibitem{biology2}
%
\begin{binproceedings}[author]
\bauthor{\bsnm{Sch{\"a}fer},~\bfnm{J.}\binits{J.}} \AND
\bauthor{\bsnm{Strimmer},~\bfnm{K.}\binits{K.}}
(\byear{2005}).
\btitle{Learning large-scale graphical {G}aussian models from genomic data}.
In \bbooktitle{Science of Complex Networks: From Biology to the
Internet and
WWW}.
\bpublisher{The American Institute of Physics}, \baddress{College Park, MD}.
\bptok{imsref}%
\end{binproceedings}
%
\endbibitem

\bibitem{ourpaper}
%
\begin{barticle}[mr]
\bauthor{\bsnm{Sturmfels},~\bfnm{Bernd}\binits{B.}} \AND
\bauthor{\bsnm{Uhler},~\bfnm{Caroline}\binits{C.}}
(\byear{2010}).
\btitle{Multivariate {G}aussian, semidefinite matrix completion, and convex
algebraic geometry}.
\bjournal{Ann. Inst. Statist. Math.}
\bvolume{62}
\bpages{603--638}.
\bid{doi={10.1007/s10463-010-0295-4}, issn={0020-3157}, mr={2652308}}
\bptok{imsref}%
\end{barticle}
%
\endbibitem

\bibitem{Whittaker}
%
\begin{bbook}[mr]
\bauthor{\bsnm{Whittaker},~\bfnm{Joe}\binits{J.}}
(\byear{1990}).
\btitle{Graphical Models in Applied Multivariate Statistics}.
\bpublisher{Wiley}, \baddress{Chichester}.
\bid{mr={1112133}}
\bptok{imsref}%
\end{bbook}
%
\endbibitem

\bibitem{biology1}
%
\begin{barticle}[author]
\bauthor{\bsnm{Wu},~\bfnm{X.}\binits{X.}},
\bauthor{\bsnm{Ye},~\bfnm{Y.}\binits{Y.}} \AND
\bauthor{\bsnm{Subramanian},~\bfnm{K.~R.}\binits{K.~R.}}
(\byear{2003}).
\btitle{Interactive analysis of gene interactions using graphical {G}aussian
model}.
\bjournal{ACM SIGKDD Workshop on Data Mining in Bioinformatics}
\bvolume{3}
\bpages{63--69}.
\bptok{imsref}%
\end{barticle}
%
\endbibitem

\end{thebibliography}
\end{document}